\RequirePackage{fix-cm}
\documentclass[smallextended]{svjour3}       
\smartqed  
\usepackage{graphicx}
\usepackage{multirow}
\usepackage{amsmath}
\usepackage{amssymb}
\usepackage{mathtools}

\def\Argmax{\mathop{\hbox{\rm  Argmax}}}

\def\Argmin{\mathop{\hbox{\rm  Argmin}}}

\def\bB{{\mathbf{B}}}

\def\beq{\begin{equation}}

\def\blacksquare{\hbox{\vrule width 4pt height 4pt depth 0pt}}

\def\eeq{\end{equation}}

\def\epscert{\varepsilon_{\hbox{\rm\scriptsize cert}}}
\def\epsiloncert{\varepsilon_{\hbox{\rm\scriptsize cert}}}

\def\epsilonsaddle{{\varepsilon_{\hbox{\scriptsize\rm sad}}}}

\def\fnote#1{\footnote}

\def\inter{\mathop{\hbox{\rm  int}\,}}
\def\inter{{\mathop{\hbox{\rm int}}}}

\def\maxl{\mathop {\max }\limits}

\def\minl{\mathop {\min }\limits}

\def\mypict3{\epsfxsize=220pt\epsfysize=80pt\epsfbox}

\def\Opt{{\hbox{\rm  Opt}}}

\def\SadSet{{\hbox{\rm SadSet}}}
\def\SadVal{{\hbox{\rm SadVal}}}

\def\bR{{\mathbf{R}}}

\journalname{myjournal}
\begin{document}

\title{Decomposition of Large Scale Linear Programming Problems Containing both Linking Variables and Constraints via Accuracy Certificates
}

\titlerunning{Decomposition via Accuracy Certificates} 

\author{Bruce A. Cox}


\institute{B. Cox \at
          Air Force Institute of Technology \\
              \email{BruceACox1@gmail.com}          
}

\date{3 March 2020}

\maketitle
\begin{abstract}
Several well known large scale linear programming decomposition methodologies exist.  Benders Decomposition, which covers the case where some small subset of variables link the otherwise separable subproblems.  Dantzig-Wolfe decomposition and Lagrangian decompositions, which cover the case where some few constraints link the otherwise separable subproblems, and finally the "Cross-Decomposition" originating from TJ Van Roy which enables one to deal with both linking constraints and linking variables by essentially alternating iteratively between the Benders and the Lagrangian Decomposition. In this paper we present a novel alternative to Cross-decomposition that deals with both linking constraints and linking variables through the application of accuracy certificates for black-box, sub-gradient based algorithms such as NERML.      
\end{abstract}

\section{Introduction}
\label{intro}
Our paper is motivated by the consideration of a large scale Linear Programming problem such as:
\begin{equation}\label{LP}
\text{Opt}=\underset{x=[x_1;x_2]}{min}\left\{c_1^Tx_1+c_2^Tx_2:\begin{array}{l}
A^{11} x_1+A^{12}x_2\leq b_1\\
A^{21}x_1+A^{22}x_2\leq b_2\\
\|x\|_\infty\leq R\\
\end{array}\right\}
\end{equation}

With $n = n_1 + n_2$ variables and $m = m_1 + m_2$ linear inequality constraints, and the sizes of $A^{11}$, $A^{12}$, $A^{21}$, and $A^{22}$ are $m_1 \text{ x }n_1$, $m_1 \text{ x } n_2$, $m_2 \text{ x } n_1$, and $m_2 \text{ x } n_2$ respectively.  Note that the bounds on variables $\|x\|_\infty\leq R$ are added for technical reasons and have no real practical significance since $R$ can be arbitrarily large.

In addition we make the following critical assumptions regarding the nature of the underlying problem.
\begin{itemize}
\item[(1)] It is relatively easy to solve Linear Programming problems of the form
\begin{equation}\label{oftheform}
\min\limits_{x_1}\left\{c^Tx_1:A^{11}x_1\leq b,\|x_1\|_\infty\leq R\right\}
\end{equation}

The classic, though by no means only situation, which meets this assumption is where $A^{11}$ is a block-diagonal matrix with a large number of relatively small diagonal blocks.  While in fact, the true general criteria is that both the Primal and the Dual linear programs can be solved quickly.

\item[(2)] The number of {\sl linking constraints} $m_2$ is $\ll$ the total number of constraints $m$ in (\ref{LP}), and the number of {\sl linking variables} $n_2$ is $\ll$ the total number of variables $n$ in (\ref{LP}).
\end{itemize}

For example the following graphic outlines one such case that meets both assumptions.

\includegraphics[scale=.5]{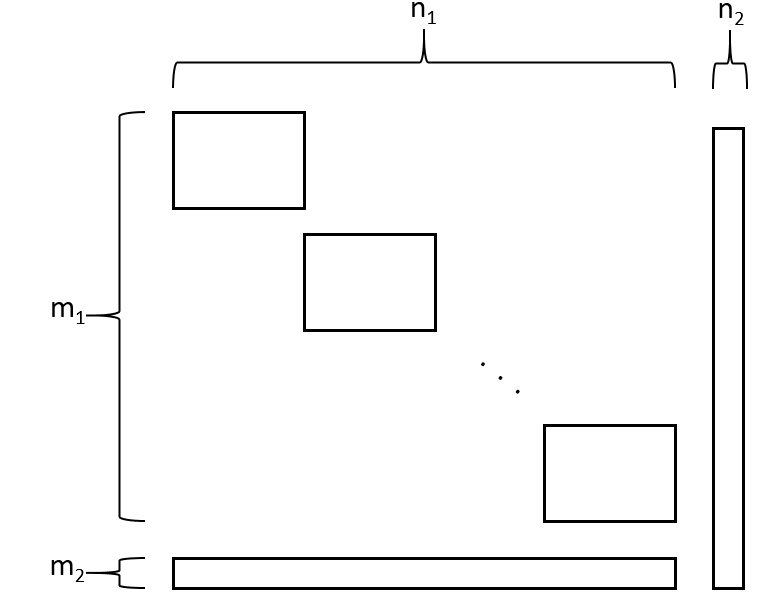}

The fundamental question is: \emph{How can we exploit the specific structure of (\ref{LP}) in order to accelerate its solution?}  This is the question we address in this paper. The organization of the paper is as follows: In section \ref{twowell} we address two well known cases of the outlined situation where we already know how to act --- one where there are no linking constraints, and one where there are no linking variables. In section \ref{sectgeneral} we address the more interesting general case where both linking variables and linking constraints are present, via a novel approach utilizing accuracy certificates.
\section{Two Well Known Cases}\label{twowell}

\subsection{Case A: No Linking Constraints}
Consider the case where there are no linking constraints, that is, $m_2=0$.
In this case, we can define the convex function
$$
f(x_2)=\min\limits_{x_1}\left\{c_1^Tx_1+c_2^Tx_2: A^{11}x_1+A^{12}x_2\leq b_1,\|x_1\|_\infty\leq R\right\}
$$
and replace the problem of interest with an equivalent problem
\begin{equation}\label{P1}
\min\limits_{X_2:\|x_2\|_\infty\leq R}f(x_2).
\end{equation}
According to assumption (1), it is relatively easy to check, given $x_2$, whether $x_2\in$ Dom $f$, and if so, to compute $f(x_2)$ and a subgradient $f'(x_2)$ of $f_2$ at $x_2$. Indeed, to this end one should solve the LP program
\begin{equation}\label{P1a}
\min\limits_{x_1}\left\{c_1^Tx_1: A^{11}x_1\leq b_1-A^{12}x_2,\|x_1\|_\infty\leq R\right\},
\end{equation}
which by our assumption is easy. If this problem is solvable, then $x_2\in$ Dom $f$, and the optimal solution $x_1(x_2)$ to the LP \eqref{P1a}, along with an optimal solution $y_1=y_1(x_2)$ to its associated dual problem
$$
\max\limits_{y_1}\left\{[A^{12}x_2-b_1]^Ty_1:[A^{11}]^Ty_1+c_1=0,y_1\geq0\right\}
$$
allows us to compute $f(x_2)$, and $f'(x_2)$ according to
$$
f(x_2)=c_1^Tx_1(x_2)+c_2^Tx_2=y_1^T[A^{12}x_2-b_1]+c_2^Tx_2,\,\,f'(x_2)=[A^{12}]^Ty_1+c_2.
$$
Conversely if the problem (\ref{P1a}) is unsolvable (which, since the feasible set of the problem clearly is bounded, can happen only if the problem \eqref{P1a} is infeasible),
we will receive, via the General Theorem on Alternative, an infeasibility certificate $y_1=y_1(x_2)$ such that
$$
y_1\geq0\ \&\ y_1^T[A^{12}x_2-b_1]>\|[A^{11}]^Ty_1\|_\infty;
$$
in this case, the linear form $\langle [A^{12}]^Ty_1,\cdot\rangle$ is nonzero and separates $x_2$ from the feasible domain $X_2$ of (\ref{P1}). Thus, (\ref{P1}) is naturally equipped with both a Separation oracle, and a First Order oracle, and as such can be solved by a black box oriented method for convex minimization, see, e.g. \cite{NERML}. By assumption, the dimension of this problem $n_2$ is $\ll$ the dimension of (\ref{LP}), which can make this outlined approach highly preferable to a direct attack on (\ref{LP}) with an LP solver.

Note that after a feasible $\epsilon$-optimal solution $x_2$ to (\ref{P1}) is found,  we automatically have at our disposal an $\epsilon$-optimal feasible solution $[x_1(x_2);x_2]$ to the problem of interest.

Note that one of the most popular implementations of this scheme is the classical {\sl Benders decomposition}, see, e.g., \cite{BerT}.

\subsection{Case B: No Linking Variables}\label{Lagrangean}
Now consider the case where there are no linking variables, that is,  $n_2=0$.
One well-known way to utilize the resulting structure is to use the well-known Dantzig-Wolfe decomposition which is intrinsically linked to the revised Primal Simplex Method, see, e.g., \cite{BerT}. Another standard way to handle
 (\ref{LP}) with no linking variables is to use {\sl Lagrangian decomposition},
that is, to associate with (\ref{LP}) its partial dual problem (where we \emph{dualize} the linking constraints) -- which we write down in the minimization form
\begin{equation}\label{P2}
\min\limits_{y_2\geq0}f(y_2), \,\,f(y_2)=\max\limits_{x_1}\left\{-[c_1+[A^{21}]^Ty_2]^Tx_1:A^{11}x_1\leq b_1,\|x_1\|_\infty\leq R\right\}+b_2^Ty_2.
\end{equation}
Observe that $f(y_2)$ is an everywhere finite convex function which can be easily equipped with a First Order oracle, since by our assumption the parametric LP problem specifying  $f$ is relatively easy to solve. Indeed the LP associated with \eqref{P2} can be rewritten as $\min \{d^Tx_1:A^{11}x_1\le b_1,\|x_1\|_1\le R\}$ where $d=c_1+[A^{2}]^Ty_2$. Denoting by $x_1(y_2)$ an optimal solution to the optimization
problem
\begin{equation}\label{probx1}
\max\limits_{x_1}\left\{-[c_1+[A^{21}]^Ty_2]^Tx_1:A^{11}x_1\leq b_1,\|x_1\|_\infty\leq R\right\},
\end{equation}
we have
\begin{equation}\label{eq447}
f(y_2)=-[c_1+[A^{21}]^Ty_2]^Tx_1(y_2)+b_2^Ty_2,\,f'(y_2)=b_2-A^{21}x_1(y_2).
\end{equation}
In many cases we can find a finite upper bound $L$ on, say, the $\|\cdot\|_1$-norm of an optimal solution to (\ref{P2}) and thus reduce this problem to a convex minimization on a simple solid, which we can solve by a black box oriented method (see e.g. \cite{NERML}). Since the dimension $m_2$ of (\ref{P2}) is small when compared to the sizes of the problem of interest (\ref{LP}), this approach can be  more attractive than a direct attack on (\ref{LP}) with an LP solver.
\par
Note that the Lagrange decomposition approach requires us to solve the nontrivial problem of recovering a good approximate solution to the problem of interest (\ref{LP}) from a good approximate solution to the problem
\begin{equation}\label{P2a}
-\Opt(L):=\min\limits_{y_2\in\bB}f(y_2),\,\,\bB=\bB_L=\left\{y_2:y_2\geq0, \sum_{i=1}^{m_2}[y_2]_i\leq L\right\}
\end{equation}
(which is nothing but the problem (\ref{P2}) with an added bound on variables). Observe that by Weak Duality $-\Opt(L)\leq \Opt$ for all $L$. One of the ways to resolve this problem is to use   accuracy certificates; to this end, we have proved the following result:

\begin{proposition}\label{lagrange}
Assume that $n_2=0$ and that {\rm (\ref{P2a})} is solved by an algorithm with accuracy certificates. By construction of the First Order oracle for the objective of this problem, at every productive step $t$ of this algorithm, the search point being $y_2^t\in \bB_L$, we have at our disposal a point
$$
x_1^t= x_1(y_2^t)\in\Argmax_{x_1}\left\{-[c_1+[A^{21}]^Ty_2^t]^Tx_1:A^{11}x_1\leq b_1,\|x_1\|_\infty\leq R\right\}.
$$
Now let $\tau$ be a step such that the accuracy certificate $\xi^\tau$ associated with the corresponding execution protocol $P_\tau=\left\{\{y_2^t,e_t\}_{t=1}^\tau, I_\tau, J_\tau\right)$ is well defined. Setting
$$
\widehat{x}_1^\tau=\sum_{t\in I_\tau}\xi^\tau_t x_1^t,
$$
we get an approximate solution to the problem of interest {\rm (\ref{LP})} such that
\begin{equation}\label{suchthat}
\begin{array}{l}
A^{11}\widehat{x}_1^\tau \leq b_1\ \&\ \|x_1^\tau\|_\infty\leq R,\\
A^{21}\widehat{x}_1^\tau\leq b_2+L^{-1}[\Opt+\epscert(\xi^\tau|P_\tau,\bB_L)+R\|c_1\|_1]{\mathbf{1}}\\
c_1^T\widehat{x}_1^\tau\leq\Opt+\epscert(\xi^\tau|P_\tau,\bB_L),\\
\end{array}
\end{equation}
where $\mathbf{1}$ is the all-ones vector of dimension $m_2$.

\par In addition, let $\widetilde{y}=[\widetilde{y}_1;\widetilde{y}_2;\widetilde{y}_+;\widetilde{y}_-]\geq0$ be the vector of optimal Lagrange multipliers for {\rm (\ref{LP})}, so that
\begin{equation}\label{that}
c_1^Tx_1+\widetilde{y}_1^T[A^{11}x_1-b_1]+\widetilde{y}_2^T[A^{21}x_1-b_2]+\widetilde{y}_+^T[x_1-R{\mathbf{1}}]+\widetilde{y}_-^T[-x_1 -R{\mathbf{1}}]\equiv \Opt\,\forall x_1.
\end{equation}
When $\ell:=L-\sum_{i=1}^{m_2}[\widetilde{y}_2]_i>0$, then, in addition to the second relation in {\rm (\ref{suchthat})}, we have
\begin{equation}\label{suchthat11}
A^{21}\widehat{x}_1^\tau\leq b_2+\ell^{-1}\epscert(\xi^\tau|P_\tau,\bB_L){\mathbf{1}}
\end{equation}
\end{proposition}

{\bf Proof.}
Let $\widehat{x}_1^\tau$ be as stated in Proposition.  Then
\begin{enumerate}
\item[i.] Clearly $A^{11}\widehat{x}_1^\tau \leq b_1\ \&\ \|\widehat{x}_1^\tau\|_\infty\leq R$ holds since by construction $\widehat{x}_1^\tau$ is a convex combination of feasible solutions to (\ref{probx1}).
\item[ii.] We have for every $y_2\in \bB_L$:
$$
\begin{array}{l}
\epscert(\xi^\tau|P_\tau,\bB_L)\geq \sum_{t=1}^\tau\xi_t\langle e_t,y_2^t-y_2\rangle \hbox{\ [definition of $\epscert$]}\\
\geq\sum_{t\in I_\tau}\xi_t\langle f'(y_2^t),y_2^t-y_2\rangle \hbox{\ [since $\langle e_t,y_2^t-y_2\rangle\geq0$ for $t\not\in I_\tau$]}\\
=\sum_{t\in I_\tau}\xi_t\langle b_2-A^{21}x_1^t,y_2^t-y_2\rangle \hbox{\ [see (\ref{eq447})]}\\
=\sum_{t\in I_\tau}\xi_t\left[\langle b_2-A^{21}x_1^t,y_2^t\rangle-\langle b_2-A^{21}x_1^t,y_2\rangle\right]\\
=\sum_{t\in I_\tau}\xi_t\left[f(y_2^t)+c_1^Tx_1^t-\langle b_2-A^{21}x_1^t,y_2\rangle\right]\hbox{\ [see (\ref{eq447})]}\\
=\sum_{t\in I_\tau}\xi_tf(y_2^t)+c_1^T\widehat{x}_1^\tau -\langle b_2-A^{21}\widehat{x}_1^\tau,y_2\rangle\\
\end{array}
$$
whence, taking into account that $f(y_2^t)\geq -\Opt(L)$ for $t\in I_\tau$,
$$
\langle A^{21}\widehat{x}_1^\tau-b_2,y_2\rangle\leq \Opt(L)+\epscert(\xi^\tau|P_\tau,\bB_L)-c_1^T\widehat{x}_1^\tau. \eqno{(!)}
$$
Relation (!) holds true for all $y_2\in\bB_L$; recalling what $\bB_L$ is and maximizing the left hand side in (!) over $y_2\in\bB_L$ with $\sum_i(y_2)_i=L$, we get
$$
\begin{array}{l}
L\max_i[A^{21}\widehat{x}_1^\tau-b_2]_i\leq \Opt(L)+\epscert(\xi^\tau|P_\tau,\bB_L)-c_1^T\widehat{x}_1^\tau\\
\leq \Opt(L)+\epscert(\xi^\tau|P_\tau,\bB_L)-\|c_1\|_1R\\
\end{array}
$$
which is nothing but the second relation in (\ref{suchthat}). Setting in (!) $y_2=0$, we get
$$
c_1^T\widehat{x}_1^\tau\leq \Opt(L)+\epscert(\xi^\tau|P_\tau,\bB_L).
$$
Taking into account that by weak duality $\Opt\geq \sup\limits_{y_2\geq0}[-f(y_2)]\geq \sup\limits_{y_2\in\bB_L}[-f(y_2)]=\Opt(L)$, we arrive at the third relation in (\ref{suchthat}).
\item[iii.]  To prove (\ref{suchthat11}), we can assume w.l.o.g. that the set $K=\{i:[A^{21}\widehat{x}_1^\tau-b_2]_i>0\}$ is nonempty, since otherwise (\ref{suchthat11}) is evident.
Setting in (\ref{that}) $x_1=\widehat{x}_1^\tau$ and taking into account the first relation in (\ref{suchthat}), we get
$$
c_1^T\widehat{x}_1^\tau+\widetilde{y}_2^T[A^{21}\widehat{x}_1^\tau-b_2]\geq\Opt\geq\Opt(L),
$$
so that
$$
\Opt(L)-c_1^T\widehat{x}_1^\tau\leq \widetilde{y}_2^T[A^{21}\widehat{x}_1^\tau-b_2],
$$
which combines with (!) to imply that for all $y_2\in\bB_L$ one has
$$
\begin{array}{l}
\langle y_2,A^{21}\widehat{x}_1^\tau-b_2\rangle\leq \epscert(\xi^\tau|P_\tau,\bB_L)+ \widetilde{y}_2^T[A^{21}\widehat{x}_1^\tau-b_2]\\
\leq \epscert(\xi^\tau|P_\tau,\bB_L)+\left[\sum_{i=1}^{m_2}[\widetilde{y}_2]_i\right]\max_i[A^{12}\widehat{x}_1^\tau-b_2]_i.\\
\end{array}
$$
Let $\mu:=\max_i[A^{12}\widehat{x}_1^\tau-b_2]_i=[A^{12}\widehat{x}_1^\tau-b_2]_{i_*}$, so that the above relation reads
$$
\langle y_2,b_2-A^{21}\widehat{x}_1^\tau\rangle\leq \epscert(\xi^\tau|P_\tau,\bB_L)+ \mu \left[\sum_{i=1}^{m_2}[\widetilde{y}_2]_i\right].
$$
Setting the $i_*$-th coordinate of $y_2$ to the value $L$, and the remaining coordinates  -- to the value 0, we get from the latter inequality
$$
\mu L\leq \epscert(\xi^\tau|P_\tau,\bB_L)+ \mu \left[\sum_{i=1}^{m_2}[\widetilde{y}_2]_i\right],
$$
and (\ref{suchthat11}) follows. \blacksquare
\end{enumerate}
\section{General Case}\label{sectgeneral}
Now consider the general case, when there are both linking variables and linking constraints. To the best of our knowledge,
the only decomposition scheme proposed for this case is ``cross decomposition'' originating from T.J. Van Roy \cite{VanRoy,VanRoy1}, see also
\cite{ChBY,Hol90,Hol92,Hol94,Hol97,Brick} and references therein. In this scheme, essentially, one alternates iteratively between the Benders and the Lagrange decompositions. To the best of our understanding, no complexity results for this scheme are known.	\footnote{This being said, note that the primary motivation and application of cross decomposition is the Mixed Integer version of (\ref{LP}), where the linking variables $x_2$ are subject to additional constraints of integrality.}
\subsection{Assumptions and Approach}
Consider the Lagrange function of problem \eqref{LP}, assuming that both linking constraints and linking variables exist.  This bilinear function is:
\begin{equation}\label{bilinear}
\Psi(x_1,x_2,y_1,y_2)=
\langle c_1,x_1\rangle+\langle c_2,x_2\rangle +\langle y_1,A^{11}x_1+A^{12}x_2-b_1\rangle + \langle y_2,A^{21}x_1+A^{22}x_2-b_2\rangle
\end{equation}
This function leads us naturally to consider the bilinear convex-concave saddle point problem equivalent to problem (\ref{LP}).
\begin{equation}\label{Psi1}
\inf\limits_{x=[x_1;x_2]\in X_1\times X_2}\sup\limits_{y=[y_1;y_2]\geq0}\Psi(x_1,x_2,y_1,y_2),\\
\end{equation}
where $X_1=\{x_1:\|x_1\|_\infty\leq R\}$, and $X_2=\{x_2:\|x_2\|_\infty\leq R\}$ are simple solids. Assume, as is usually the case, that we can
shrink the $y$-domain of this saddle point problem, which is a nonnegative orthant, to a direct product $Y=Y_1\times Y_2$ of two simple solids  which are ``large enough'' to ensure that the saddle point problem
\begin{equation}\label{Psi2}
\min\limits_{x=[x_1;x_2]\in X_1\times X_2}\max\limits_{y=[y_1;y_2]\in Y_1\times Y_2} \Psi(x_1,x_2,y_1,y_2)
\end{equation}
is equivalent to (\ref{Psi1}) and thus -- to the LP of interest (or approximates (\ref{Psi1}) ``well enough'').
\paragraph{Assumption.} At the beginning of this chapter we have assumed that it is easy
to solve the LPs of the form $\min\limits_x\{c^Tx_1:A^{11}x_1\leq b,\|x\|_\infty\leq R\}$. From now on we modify this assumption, specifically, assume that
 \begin{quote}
 {\bf A}: {\sl Given $x_2\in X_2,y_2\in Y_2$, it is easy to solve the saddle point problem}
 \begin{equation}\label{SPeasy}
 \begin{array}{c}
 \min\limits_{x_1\in X_1}\max\limits_{y_1\in Y_1}\left\{\Psi(x_1,x_2,y_1,y_2)=\langle c+[A^{11}]^Ty_1,x_1\rangle -\langle b,y_1\rangle\right\},\\
 \left[c=c_1+[A^{12}]^Ty_2,\,b=b_1-A^{12}x_2\right]
 \end{array}
 \end{equation}
 \end{quote}
 Note that the original assumption is nothing but the modified one with $Y_1=\bR^{m_1}_+$. For our machinery to work, we need $Y_1$ to be bounded, this is where the modification comes from. Note also that in typical situations where the unmodified assumption holds true,
 so is the modified one, provided that $Y_1$ is chosen properly. For example, this is so in the case considered in section \ref{chap4motiv}, where the problem $\min_{x_1\in X_1}\{c^Tx_1:A^{11}x_1\leq b,\|x_1\|_\infty\leq R\}$ is just a collection of a large number of independent of each other LPs, or, which is the same, $x_1=[x_{11};...;x_{1K}]$ is a collection of a large number $K$ of  components $x_{1k}$ with dimensions low as compared to $n_1$, and $A^{11}$ is block-diagonal: $A^{11}x_1=[A^{11}_1x_{11};A^{11}_2x_{12};...;A^{11}_Kx_{1K}]$, with the dimensions of the blocks $A^{11}_kx_{1k}$ low as compared to $m_1$. In this case, choosing $Y_1$ as a large box: $$
 Y_1=\{y=[y_{11};...;y_{1K}]: 0\leq y_{1k}\leq R_k[1;...;1]\}$$
 with the blocks $y_{1k}$ in $y_1$ corresponding to the blocks $A^{11}_kx_{1k}$ in $A^{11}x_1$, the saddle point problem in {\bf A} decomposes into $K$ independent low dimensional bilinear saddle point problems of the form
 $$
 \min\limits_{x_{1k},\|x_{1k}\|_\infty\leq R}\max\limits_{y_{1k}:0\leq y_{1k}\leq R_k[1;...;1]} \left[\langle p_k,x_{1k}\rangle
 +\langle q_k,y_{1k}\rangle +\langle y_{1k},A^{11}_kx_{1k}\rangle\right],
 $$
 thus making assumption {\bf A} quite realistic.
\paragraph{Approach.} On a closer inspection, Assumption {\bf A} says that it is easy to build the First Order oracle for the (clearly convex-concave)
function
$$
\widetilde{\Psi}(x_2,y_2)=\min\limits_{x_1\in X_1}\max\limits_{y_1\in Y_1}\Psi(x_1,x_2,y_1,y_2),
$$
so that the saddle point problem
\begin{equation}\label{reduced}
\min\limits_{x_2\in X_2}\max\limits_{y_2\in Y_2}\widetilde{\Psi}(x_2,y_2)
\end{equation}
(which, by assumption, is of dimensions much smaller than those of the problem of interest) is well suited for solving by a black-box-oriented method, like NERML or perhaps even the Ellipsoid method\footnote{The latter is quite realistic, provided that the sizes $n_2$ and $m_2$ are in the range of tens; note that there are meaningful applications where $n_2$ and $m_2$ are indeed in the range of tens, while $m_1$ and $n_1$ are in the range of many thousands.}.
\par
 When the sizes of (\ref{reduced}) are small as compared to those of (\ref{Psi2}) and elimination of $x_1,y_1$, that is, computing the first order information for $\widetilde{\Psi}$, is cheap, the outlined approach can be computationally much more attractive than a direct attack on (\ref{Psi1}) (or, equivalently, on (\ref{LP})). Note that with this approach, we in fact {\sl simultaneously eliminate both linking variables and linking constraints} (and do it just once), in a sharp contrast to cross decomposition, where we iteratively eliminate, in an alternating fashion, either linking variables, or linking constraints.
 \par
 The crucial (and nontrivial) question underlying the outlined approach is how to recover a good solution to the problem of actual interest (\ref{Psi2}) (or, which is the same, (\ref{LP})) from a good solution to the induced saddle point problem (\ref{reduced}). Our major goal in the sequel  is to demonstrate that {\sl when solving {\rm (\ref{reduced})} by an algorithm with accuracy certificates, these certificates allow to convert $\epsilon$-solutions to {\rm (\ref{reduced})} into $\epsilon$-solutions of the saddle point problem of actual interest {\rm (\ref{Psi2})} and thus -- to $\epsilon$-solutions to the original problem  {\rm (\ref{LP})}}.
\subsection{Induced Pairs of Saddle Point Problems}\label{sect:induced}
We intend to consider a situation which is an extension of the one considered in the previous section. Specifically, we intend to investigate saddle point problems {\sl induced} by larger saddle point problems. The general setting is as follows.

Let $X_i\subset \bR^{n_{x_i}}$, and $Y_i\subset\bR^{n_{y_i}}$, $i=1,2$, be solids. We associate with these solids the sets
\begin{equation}\label{sets}
\begin{array}{cc}
X=X_1\times X_2\subset\bR^{n_x}=\bR^{n_{x_1}}\times \bR^{n_{x_2}},\quad
Y=Y_1\times Y_2\subset\bR^{n_y}=\bR^{n_{y_1}}\times \bR^{n_{y_2}},\\
Z=X\times Y\subset\bR^{n_z}=\bR^{n_x}\times\bR^{n_y}.\\
\end{array}
\end{equation}
Let $\Psi(z):Z\to\bR$ be a continuous function; with slight abuse
of notation, we shall denote this function also $\Psi(x,y)$,
$\Psi(x_1,x_2,y_1,y_2)$, etc. From now on, we make the following
\begin{quote}
{\bf Assumption B.} {\sl $\Psi$ is convex in $x\in X$ and concave
in $y\in Y$.}
\end{quote}
Note that Assumption {\bf B} is automatically satisfied when $\Psi$ is bilinear in $x,y$, as it is the case in (\ref{bilinear}) -- (\ref{Psi2}).
\par
Function $\Psi$ gives rise to two functions as follows:
$$
\begin{array}{l}
\Psi_1(x_1,y_1)=\minl_{x_2\in X_2}\maxl_{y_2\in Y_2}\Psi(x_1,x_2,y_1,y_2):X_1\times Y_1\to\bR,\\
\Psi_2(x_2,y_2)=\minl_{x_1\in X_1}\maxl_{y_1\in Y_1}\Psi(x_1,x_2,y_1,y_2):X_2\times Y_2\to\bR\\
\end{array}
$$
Since $Z$ is a compact set and $\Psi$ is continuous on $Z$, these
functions are continuous on the respective solids $Z_i=X_i\times
Y_i$.

\begin{lemma}\label{lem1} $\Psi_i(x_i,y_i)$ is convex in $x_i\in
X_i$ and is concave in $y_i\in Y_i$, $i=1,2$.
\end{lemma}
{\bf Proof.} Indeed, we have, that $\Psi_1(x_1,y_1)=\min_{x_2\in
X_2}\left[\max_{y_2\in Y_2}\Psi(x_1,x_2,y_1,y_2)\right]$ is
concave in $y_1$, since the function in brackets $[\,]$ is concave
in $y_1$ due to the concavity of $\Psi(x,y)$ in $y$. The convexity
of $\Psi_1(x_1,y_1)$ in $x_1$ follows via similar argument from
the representation $\Psi_1(x_1,y_1)=\max_{y_2\in
Y_2}\left[\min_{x_2\in X_2}\Psi(x_1,x_2,y_1,y_2)\right]$, where
the interchange of $\min$ and $\max$ is legitimate due to the fact
that $\Psi(x_1,x_2,y_1,y_2)$ is convex-concave and continuous in
$x_2,y_2$ and $X_2,Y_2$ are convex compact sets (from now on, we
skip similar justifications of swapping the order of $\min$ and
$\max$). The proof for $\Psi_2(x_2,y_2)$ follows a symmetric line of arguments. \blacksquare
\par
We are about to consider two convex-concave saddle point problems
\begin{equation}\label{PD}
\begin{array}{lc}
(P): \SadVal_{X_1\times Y_1}(\Psi_1)=\minl_{x_1\in X_1}\maxl_{y_1\in Y_1}\Psi_1(x_1,y_1)\\
(D): \SadVal_{X_2\times Y_2}(\Psi_2)=\minl_{x_2\in X_2}\maxl_{y_2\in Y_2}\Psi_2(x_2,y_2)\\
\end{array}
\end{equation}
which we call the problems {\sl induced} by the master saddle
point problem
\begin{equation}\label{master}
(M): \SadVal_{X\times Y}(\Psi)=\min_{x=(x_1,x_2)\in X=X_1\times
X_2}\max_{y=(y_1,y_2)\in Y=Y_1\times Y_2}\Psi(x,y).
\end{equation}
We start with the following observation:

\begin{proposition}\label{lem2} One has
\begin{equation}\label{onehas1}
\SadVal_{X_1\times Y_1}(\Psi_1)=\SadVal_{X_2\times
Y_2}(\Psi_2)=\SadVal_{X\times Y}(\Psi)
\end{equation}
and
\begin{equation}\label{onehas2}
\SadSet_{X\times Y}(\Psi)\subset\SadSet_{X_1\times Y_1}(\Psi_1)\times \SadSet_{X_2\times Y_2}(\Psi_2),
\end{equation}
where $\SadSet_{U\times V}(\Psi)$ is the set of saddle points of the function
$\Psi(\cdot)$ on $U\times V$.
\par
Moreover, for every $(x_1,x_2,y_1,y_2)\in Z$ one has
\begin{equation}\label{onehas3}
\epsilonsaddle((x_i,y_i)|\Psi_i,X_i,Y_i)\leq\epsilonsaddle((x,y)|\Psi,X,Y),\,\,i=1,2.
\end{equation}
\end{proposition}

{\bf Proof.}
Let $(x^*,y^*)$ be a saddle point of $\Psi$ on $X\times Y$. We have for $\Psi_1$, that
$$
\begin{array}{l}
\Psi_1(x_1^*,y_1)=\min_{x_2\in X_2}\max_{y_2\in Y_2}
\Psi(x_1^*,x_2,y_1,y_2) =\max_{y_2\in Y_2}\min_{x_2\in
X_2}\Psi(x_1^*,x_2,y_1,y_2)\\
\leq\max_{y_2\in Y_2}\Psi(x_1^*,x_2^*,y_1,y_2)\leq
\Psi(x_1^*,x_2^*,y_1^*,y_2^*),\\
\end{array}
$$
and
$$
\Psi_1(x_1,y_1^*)=\min_{x_2\in X_2}\max_{y_2\in Y_2}
\Psi(x_1,x_2,y_1^*,y_2)\geq\min_{x_2\in X_2}
\Psi(x_1,x_2,y_1^*,y_2^*)\geq \Psi(x_1^*,x_2^*,y_1^*,y_2^*),
$$
so that
\begin{equation}\label{foralll}
\forall (x_1\in X_1,y_1\in Y_1):
\Psi_1(x_1^*,y_1)\leq\Psi(x^*,y^*)\leq \Psi_1(x_1,y_1^*).
\end{equation}
By the standard definition of saddle points as applied to $\Psi_1$ it follows that $\Psi_1(x_1^*,y_1^*)=\Psi(x^*,y^*)$,
(and thus that $\SadVal(\Psi_1)=\SadVal(\Psi)$) which combines with (\ref{foralll}) to imply that $(x_1^*,y_1^*)\in\SadSet(\Psi_1)$.
By ``symmetric'' reasoning, $(x_2^*,y_2^*)\in\SadSet(\Psi_2)$ and $\SadVal(\Psi_2)=\SadVal(\Psi)$.  Thus (\ref{onehas1}) and (\ref{onehas2}) are proved.

\par Now let $(x_1,x_2,y_1,y_2)\in Z$. We have by \eqref{PrimalDualPair}
$$
\begin{array}{l}
\epsilonsaddle((x_1,y_1)|\Psi_1,X_1,Y_1)=\max_{\eta_1\in
Y_1}\Psi_1(x_1,\eta_1)-\min_{\xi_1\in X_1}\Psi_1(\xi_1,y_1)\\
=\max_{\eta_1\in Y_1}\min_{\xi_2\in X_2}\max_{\eta_2\in
Y_2}\Psi(x_1,\xi_2,\eta_1,\eta_2)\\
\quad-\min_{\xi_1\in
X_1}\min_{\xi_2\in
X_2}\max_{\eta_2\in Y_2}\Psi(\xi_1,\xi_2,y_1,\eta_2)\\
\leq \max_{\eta_1\in Y_1}\max_{\eta_2\in
Y_2}\Psi(x_1,x_2,\eta_1,\eta_2)\\
\quad-\min_{\xi_1\in X_1}\min_{\xi_2\in
X_2}\Psi(\xi_1,\xi_2,y_1,y_2)\\
=\epsilonsaddle((x,y)|\Psi,X,Y),\\
\end{array}
$$
and similarly for $\Psi_2$. (\ref{onehas3}) is proved.
\blacksquare
\subsection{Recovering Approximate Solutions to the Master Problem: Goal and Assumptions}
\paragraph{The goal.} Consider a master saddle point problem (\ref{master}) along with the induced problems (\ref{PD}). By Proposition \ref{lem2}, (specifically by \eqref{onehas3}) we can easily extract good approximate solutions to each of the induced problems from a good approximate solution to the master problem. The question is: To what extent is the opposite true?

Specifically, assume we have at our disposal a first-order method capable of solving to within a desired accuracy one of the induced problems, say problem $(P)$, and our goal is to extract from this
solution a good approximate solution to the master problem (and thus, by Proposition \ref{lem2}, to problem $(D)$ as well). When and how could we achieve this goal?

We will demonstrate that this goal is achievable, provided that \begin{itemize}
\item the first order information used by the algorithm in question satisfies some not too restrictive
technical assumptions, and that
\item we have at our disposal, not only the approximate solution to $(P)$ to be converted into
approximate solutions to $(M)$ and $(D)$, but also an accuracy certificate for this solution.
\end{itemize}

\paragraph{Preliminaries.} We start with some technical issues.   From now on, we make
\begin{quote}
{\bf Assumption C.}
{\sl $\Psi(x_1,x_2,y_1,y_2)$ is not only continuous convex-concave on $Z=X\times Y$, but
\begin{itemize}
\item $\Psi$ is differentiable in $x_2\in X_2$ whenever $x_1\in\inter X_1$, $y_1\in\inter Y_1$ and $y_2\in Y_2$, the derivative being continuous in $x\in \inter X_1\times X_2$ for every $y\in\inter
Y_1\times Y_2$;
\item $\Psi$ is differentiable in $y_2\in Y_2$ whenever $x_1\in\inter X_1$, $y_1\in\inter Y_1$ and $x_2\in X_2$, the derivative being continuous in $y\in \inter Y_1\times Y_2$ for every $x\in\inter X_1\times X_2$.
\end{itemize}}

\end{quote}
Under Assumption {\bf C}, for every point $(x=(x_1,x_2),y=(y_1,y_2))\in Z$ with $(x_1,y_1)\in \inter
X_1\times \inter Y_1$, function $\Psi(\cdot,y)$ at the point $x$ admits a {\sl regular} subgradient $\Psi^\prime_x(x,y)$ -- that is a subgradient whose projection onto the subspace $\bR^{n_{x_2}}$
equals the gradient of the continuously differentiable function $\Psi(x_1,\cdot,y)$ at the point $x_2$.
\begin{quote}
{\small Indeed, let $x_2^t\in\inter X_2$ be such that $x_2^t\to x_2$ as $t\to\infty$, and let $g^t$ be a
subgradient of $\Psi(\cdot,y)$ at the point $x^t=(x_1,x_2^t)\in\inter X$. By evident reasons, such a
subgradient is automatically regular.  Besides this, the vectors $g^t$ form a bounded sequence, since $\Psi(\xi_1,\xi_2,y)$ is uniformly in $\xi_1\in V$ Lipschitz continuous in $\xi_2\in X_2$, $V$ being a neighbourhood of $x_1$ with the closure belonging to $\inter X_1$,  and is uniformly in $\xi_2\in X_2$ Lipschitz continuous in $\xi_1\in V$. Passing to a subsequence. we may assume that $g^t$ has a limit $g$ as $t\to\infty$; by construction, $g$ is a subgradient of $\Psi(\cdot,y)$ at the point
$x$, and its projection on $\bR^{n_{x_2}}$ is $\lim_{t\to\infty}\nabla_s\big|_{s=x_2^t}\Psi(x_1,s,y)=\nabla_s\big|_{s=x_2}\Psi(x_1,s,y)$,
as required from a regular subgradient.}
\end{quote}

We define similarly the notion of a {\sl regular} supergradient $\Psi^\prime_y(x,y)$ in $y$ (a supergradient of the concave function $\Psi(x,\cdot)$ at the point $y$ such that the projection of this supergradient onto $\bR^{n_{y_2}}$ is the gradient of $\Psi(x,y_1,\cdot)$ at the point $y_2$). Such a
supergradient also exists, provided that $(x_1,y_1)\in\inter X_1\times \inter Y_1$, and $(x_2,y_2)\in X_2\times Y_2$.
\begin{remark}\label{normality} Note that Assumption {\bf C} is automatically satisfied when $\Psi$ is convex-concave and continuously differentiable (as it is the case, e.g., when $\Psi$ is bilinear in $x,y$, cf. {\rm (\ref{bilinear}), (\ref{Psi2})}). In this case, choosing as subgradients of  $\Psi$ w.r.t. $x$ the corresponding partial gradients, and similarly for supergradients of $\Psi$ in $y$, we automatically end up with regular sub- and supergradients.
\end{remark}
\begin{lemma}\label{lem3}
Given
$(\bar{x}_1,\bar{y}_1)\in \inter X_1\times \inter Y_1$, let
$(\bar{x}_2,\bar{y}_2)$ be a saddle point of the convex-concave
continuous function
$\Psi_{\bar{x}_1,\bar{y}_1}(x_2,y_2)=\Psi(\bar{x}_1,x_2,\bar{y}_1,y_2)$
on $X_2\times Y_2$, and let $ \Psi^\prime_x(\bar{z})$,
$\Psi^\prime_y(\bar{z})$ be regular sub- and supergradients of
$\Psi$ in $x$ and in $y$, respectively, computed at the point
$\bar{z}=(\bar{x}_1,\bar{x}_2,\bar{y}_1,\bar{y}_2)$. Let, further,
$\Psi^\prime_{1,x}(\bar{x}_1,\bar{y}_1)$ be the projection of  $
\Psi^\prime_x(\bar{z})$ onto $\bR^{n_{x_1}}$, and
$\Psi^\prime_{1,y}(\bar{x}_1,\bar{y}_1)$ be the projection of $
\Psi^\prime_y(\bar{z})$ onto $\bR^{n_{y_1}}$. Then
$\Psi^\prime_{1,x}(\bar{x}_1,\bar{y}_1)$ is a subgradient of the
convex function $\Psi_1(x_1,\bar{y}_1)$ of $x_1\in X_1$ at the
point $x_1=\bar{x}_1$, and
$\Psi^\prime_{1,y}(\bar{x}_1,\bar{y}_1)$ is a supergradient of the
concave function $\Psi_1(\bar{x}_1,y_1)$ of $y_1\in Y_1$ at the
point $y_1=\bar{y}_1$.
\end{lemma}

{\bf Proof.}
For $x_1\in X_1$, we have
$$
\begin{array}{l}
\Psi_1(x_1,\bar{y}_1)
=\min_{x_2\in X_2}\max_{y_2\in Y_2}\Psi(x_1,x_2,\bar{y}_1,y_2)
\geq \min_{x_2\in X_2}\Psi(x_1,x_2,\bar{y})\\
\geq \min_{x_2\in X_2}\left[\Psi(\bar{x},\bar{y})+ \langle \Psi^\prime_x(\bar{z}),x-\bar{x}\rangle\right]
\hbox{\ [since $\Psi(\cdot,\bar{y})$ is convex]}\\
=\underbrace{\Psi_1(\bar{x}_1,\bar{y}_1)}_{=\Psi(\bar{x},\bar{y})}+\min_{x_2\in X_2}\left[\langle
\Psi^\prime_{1,x}(\bar{x}_1,\bar{y}_1),x_1-\bar{x}_1\rangle+\langle
\nabla_s\big|_{s=\bar{x}_2}\Psi(\bar{x}_1,s,\bar{y}),x_2-\bar{x}_2\rangle\right]\\
\multicolumn{1}{r}{\hbox{[since $\Psi^\prime_x(\bar{z})$ is regular]}}\\
\geq \Psi_1(\bar{x}_1,\bar{y}_1)+\min_{x_2\in X_2}\left[\langle \Psi^\prime_{1,x}(\bar{x}_1,\bar{y}_1),x_1-\bar{x}_1\rangle\right]\\
\multicolumn{1}{r}{\left[\hbox{\footnotesize\begin{tabular}{l} since $\bar{x}_2\in\Argmin_{x_2\in X_2}f(x_2)$,
$f(x_2)=\Psi(\bar{x}_1,x_2,\bar{y}_1,\bar{y}_2)$, and $f(x_2)$ is differentiable\\ at $\bar{x}_2$, so that $\langle
\nabla_s\big|_{s=\bar{x}_2}\Psi(\bar{x}_1,s,\bar{y}),x_2-\bar{x}_2\rangle\geq0$ whenever $x_2\in X_2$\\
\end{tabular}}\right]}\\
= \Psi_1(\bar{x}_1,\bar{y}_1)+\langle \Psi^\prime_{1,x}(\bar{x}_1,\bar{y}_1),x_1-\bar{x}_1\rangle;\\
\end{array}
$$
The concluding inequality says that
$\Psi^\prime_{1,x}(\bar{x}_1,\bar{y}_1)$ indeed is a subgradient
of $\Psi_1(x_1,\bar{y}_1)$ in $x_1\in X_1$ evaluated at
$x_1=\bar{x}_1$. The ``symmetric'' reasoning proves the
``supergradient'' part of the statement. \blacksquare
\subsection{Recovering Approximate Solutions to the Master Problem: Construction and Main Result}\label{mainr}
In the above described situation, let us assume that we have access to Separation oracles for $X_1$ and $Y_1$ (and thus -- to a Separation oracle for $Z_1=X_1\times Y_1$). We also assume that we have access to a $\Phi$-oracle, where $\Phi:\inter Z_1\to\bR^{n_{x_1}}\times \bR^{n_{y_1}}$ is the monotone mapping
associated with the convex-concave saddle point problem
\begin{equation}\label{eq1}
\max_{x_1\in X_1}\min_{y_1\in Y_1}\Psi_1(x_1,y_1),
\end{equation}
specifically, as follows:

\begin{quote}
Given on input $(x_1,y_1)\in\inter X_1\times\inter Y_1$, the $\Phi$-oracle
\begin{itemize}
\item solves the saddle point problem
	\begin{equation}\label{auxproblem}
	\min_{\xi_2\in X_2}\max_{\eta_2\in Y_2}\Psi(x_1,\xi_2,y_1,\eta_2)
	\end{equation}
and computes a saddle point $(x_2,y_2)$ of this problem, along with
	\begin{itemize}
	\item the projection $e_{x_1}$ of the regular subgradient of $\Psi(\xi,y_1,y_2)$ in $\xi\in X$
	computed at the point $\xi=(x_1,x_2)$, onto the space $\bR^{n_{x_1}}$;
	\item the projection $-e_{y_1}$ of the regular supergradient of $\Psi(x_1,x_2,\eta)$ in $\eta\in Y$ 		 computed at the point $\eta=(y_1,y_2)$, onto the space $\bR^{n_{y_1}}$.
	\end{itemize}
\item returns the pair $(x_2,y_2)$ and the vector $\Phi(x_1,y_1)=(e_{x_1},e_{y_1})$, thus, by Lemma \ref{lem3}, reporting the value at $(x_1,y_1)$ of the monotone mapping associated with the saddle point problem (\ref{eq1}).
\end{itemize}
\end{quote}

Assume that we have built a $\tau$-point execution protocol $P_\tau=\{(z_1^t,e^t)\}_{t=1}^\tau$, where $z_1^t=(x_1^t,y_1^t)$ are the search points, partitioned into those which are strictly feasible ($z_1^t\in\inter Z_1\Leftrightarrow t\in I_\tau$) and all the remaining search points
($z_1^t\not\in \inter Z_1\Leftrightarrow t\in J_\tau$). Also assume that $e^t$ is either $\Phi(z_1^t)$ (this is so when $t\in I_\tau$), or $e^t$ is a nonzero separator of $z_1^t$ and $Z_1$ (this is so when $t\in J_\tau$). According to the construction of the $\Phi$-oracle, this protocol can be augmented with pairs $z_2^t=(x_2^t,y_2^t)$, $t\in I_\tau$, reported by the $\Phi$-oracle at the productive steps (those from $I_\tau$). Our
main result is as follows:

\begin{theorem}\label{the1}
Let $\bB$ be a solid containing $Z_1$, and let $P_\tau$ be an execution protocol which admits an accuracy certificate $\zeta$. Given this certificate, let us set
$$
\widehat{x}_i^\tau=\sum_{t\in
I_\tau}\zeta_tx_i^t,\,\,\widehat{y}_i^\tau=\sum_{t\in
I_\tau}\zeta_ty_i^t,\,i=1,2.
$$
Then $(\widehat{x}^\tau,\widehat{y}^\tau)\in Z$ and
\begin{equation}\label{main}
\epsilonsaddle((\widehat{x}^\tau,\widehat{y}^\tau)|\Psi,X,Y)\leq\epsiloncert(\zeta|P_\tau,\bB),
\end{equation}
whence, by Proposition \ref{lem2}, also
\begin{equation}\label{maina}
\epsilonsaddle((\widehat{x}_i^\tau,\widehat{y}_i^\tau)|\Psi_i,X_i,Y_i)\leq\epsiloncert(\zeta|P_\tau,\bB),\,\,i=1,2.
\end{equation}
\end{theorem}

{\bf Proof.}
For $t\in I_\tau$, let $z^t=(x_1^t,x_2^t,y_1^t,y_2^t)=(x^t,y^t)$. Recall that for $t\in
I_\tau$, $e^t_{x_1}$ is the projection onto $\bR^{n_{x_1}}$ of a subgradient $\Psi^\prime_x(x^t,y^t)$ of the function $\Psi(\cdot,y^t)$ computed at the point $x^t$, and the projection
$e^t_{x_2}$ of this subgradient onto $\bR^{n_{x_2}}$ is the vector $\nabla_{x_2}\big|_{x_2=x_2^t}\Psi(x_1^t,x_2,y_1^t,y_2^t)$, whence
	\begin{equation}\label{eq2}
	\langle e^t_{x_2},x_2-x_2^t\rangle\geq0\,\,\forall x_2\in X_2,
	\end{equation}
due to the fact that $\Psi(x_1^t,x_2,y_1^t,y_2^t)$ attains its minimum over $x_2\in X_2$ at the point $x_2^t$. Similarly, for $t\in I_\tau$, $e^t_{y_1}$ is the projection onto $\bR^{n_{y_1}}$
of a subgradient $-\Psi^\prime_y(x^t,y^t)$ of the function $-\Psi(x^t,\cdot)$ computed at the point $y^t$, and the projection $e^t_{y_2}$ of this subgradient onto $\bR^{n_{y_2}}$ is the vector
$-\nabla_{y_2}\big|_{y_2=y_2^t}\Psi(x_1^t,x_2^t,y_1^t,y_2)$, whence
	\begin{equation}\label{eq3}
	\langle e^t_{y_2},y_2-y_2^t\rangle\geq0\,\,\forall y_2\in Y_2.
	\end{equation}
Let $z=(x_1,x_2,y_1,y_2)=(x,y)\in Z$. We have
$$
\begin{array}{l}
-\epsiloncert(\zeta|P_\tau,\bB)\leq \sum_{t=1}^\tau\zeta_t\langle
e^t,z-z^t\rangle\\
\leq \sum_{t\in I_\tau}\zeta_t\langle e^t,z-z^t\rangle\hbox{\
[since $e^t$ separates $Z$ and $z^t$ for $t\in J_\tau$]}\\
=\sum_{t\in I_\tau}\zeta_t\left[\langle e^t_{x_1},x_1-x_1^t\rangle
+ \langle e^t_{y_1},y_1-y_1^t\rangle\right]
\\
\leq \sum_{t\in I_\tau}\zeta_t\left[\langle
\Psi^\prime_x(x^t,y^t),x-x^t\rangle + \langle
-\Psi^\prime_y(x^t,y^t),y-y^t\rangle\right] \hbox{\ [by
(\ref{eq2}), (\ref{eq3})]}\\
\leq\sum_{t\in I_\tau}\zeta_t\left[[\Psi(x,y^t)-\Psi(x^t,y^t)]
+[\Psi(x^t,y^t)-\Psi(x^t,y)]\right]\hbox{\ [since $\Psi$ is
convex-concave]}\\
=\sum_{t\in I_\tau}\zeta_t\left[\Psi(x,y^t)-\Psi(x^t,y)\right]\\
\leq\Psi(x,\widehat{y}^\tau)-\Psi(\widehat{x}^\tau,y)\hbox{\
[since $\Psi$ is convex-concave]}\\
\end{array}
$$Thus,
$$
\forall(x\in X,y\in Y):
\Psi(\widehat{x}^\tau,y)-\Psi(x,\widehat{y}^\tau)\leq\epsiloncert(\zeta|P_\tau,\bB).
$$
Taking the supremum of the left hand side in $(x,y)\in Z$, we
arrive at (\ref{main}). \blacksquare




\end{document}